\magnification=1200

\tolerance=500

\font \twelvebf=cmbx12

\hfuzz=3pt

\def \a {\alpha}

\def \Coker {{\rm Coker\,}}

\def\dem{\noindent {\bf D\'emonstration.}\enspace \nobreak }

\def \det {{\rm d\acute et\,}}

\def \dip {{\rm dp\,}}

\def \equi {\Longleftrightarrow}

\def \expl#1 {\medbreak \noindent {\bf Exemple #1.}\enspace }
\def \expls#1 {\medbreak \noindent {\bf Exemples #1.}\enspace }

\def \cExt {{{\cal E}\it xt\,}}
\def \Ext {{\rm Ext\,}}

\def \f {\varphi }
\def \fl {\rightarrow }
\def \Fl#1{\smash{\mathop{\longrightarrow}\limits^{#1}}}

\def \G {\Gamma}

\def \Hom { {\rm Hom\,}}

\def \id #1 {<\!\! #1 \!\! >}

\def \lign {\hfil \break }
\def \lf {\leftarrow} 

\def \Vf#1{\Big\downarrow \rlap{$\vcenter{\hbox{$\scriptstyle#1$}}$}}

\def \ov {\overline}

\def \rema#1 {\medbreak \noindent {\bf Remarque #1.}\enspace }
\def \remas#1 {\medbreak \noindent {\bf Remarques #1.}\enspace }

\def \rg {{ \rm rang \,}}

\def \s {\sigma}

\def \Spec {{\rm Spec\,}}

\def \T {\otimes }
\def \tarte#1 {\medbreak \noindent {\it #1.}\medbreak}

\def \titre#1{\medbreak \noindent {\bf #1.}\medbreak}

\outer \def \th #1. #2\par{ \medbreak 
\noindent {\bf#1. \enspace} {\sl#2 }\par
\ifdim \lastskip< \medskipamount \removelastskip \penalty55 \medskip \fi}

\def  \Vf#1{\Big\downarrow \rlap{$\vcenter{\hbox{$\scriptstyle#1$}}$}}
\def \vf {\downarrow}

\def \Cpont {{C^{\vphantom C}}_{\hbox {\bf .}}}
\def \Lpont {{L^{\vphantom L}}_{\hbox {\bf .}}}
\def \Mpont {{M^{\vphantom M}}_{\hbox {\bf .}}}

\def \tas {{\hbox {\bf .}}}

\def \cC {{\cal C}}
\def \cE {{\cal E}}
\def \cL {{\cal L}}
\def \cP {{\cal P}}
\def \cF {{\cal F}}

\def \cH {{\cal H}}
\def \cJ {{\cal J}}
\def \cO {{\cal O}}

\def \cM {{\cal M}}
\def \cN {{\cal N}}
\def \cQ {{\cal Q}}

\def \cS {{\cal S}}
\def \cT {{\cal T}}

\def \bP {{\bf P}}

\def \bZ {{\bf Z}}
\def \bN {{\bf N}}

\def \sN{{\sl N}}

\def \cQ {{\cal Q}}
\def \Id {{\rm Id}}

\def \sN {{\sl N}}

\def \sE {{\sl E}}

\def \ps {{\it psi}}

\def \cD {{\cal D}}

\def \cQ {{\cal Q}}

\centerline{ \twelvebf Un th\'eor\`eme de Rao pour les familles de courbes gauches}

\vskip 2 cm

\titre { Introduction}

On conna\^\i t le th\'eor\`eme de Rao pour les courbes 
 localement Cohen-Macaulay et \'equidimensionnelles  de $\bP^3$ (espace projectif
de dimension $3$  sur un corps
$k$ alg\'ebrique\-ment clos).  Rappelons que le module de Rao d'une courbe $C$  est
le  $R$-module gradu\'e de longueur finie (avec
$R = k[X,Y,Z,T]$) :
$$M_C =
\bigoplus_{n \in
\bZ} H^1 \cJ_C (n).$$  
On renvoie \`a [MDP1] ou [Mi] pour des pr\'ecisions sur ce module et sur la liaison des
courbes gauches. On a alors le th\'eor\`eme de Rao originel (cf. [R1], 1979) :

\th {Th\'eor\`eme de Rao, premi\`ere forme}. Deux courbes $C$ et $C'$ sont dans la
m\^eme classe de liaison (resp. de biliaison, ou liaison paire) si et seulement si
leurs modules de Rao sont isomorphes ou duaux, \`a d\'ecalage pr\`es (resp.
isomorphes, \`a d\'ecalage pr\`es). 

 Rao a \'enonc\'e une variante de son th\'eor\`eme (cf. [R2] 1981) que l'on peut
formuler plus ais\'ement avec la notion de r\'esolution de type {\sl N}
(resp. {\sl E}) (cf. [MDP1] II) : si
$C$ est une courbe,  une telle r\'esolution est une suite exacte de faisceaux sur
$\bP^3$:
$$0 \fl \cP \fl \cN \fl \cJ_C \fl 0 \quad {\rm (resp.} \quad 0 \fl \cE \fl \cF\fl
\cJ_C
\fl 0)$$
o\`u $\cP$ et $\cF$ sont dissoci\'es (i.e. sommes d'inversibles) et o\`u $\cN$ et
$\cE$ sont des faisceaux localement libres v\'erifiant $H^2_* \cN=0$ (resp. $H^1_*
\cE=0$).  Ces faisceaux sont bien d\'etermin\'es \`a isomorphisme stable pr\`es,
c'est-\`a-dire \`a adjonction de sommes directes d'inversibles pr\`es. On a ainsi une
application qui \`a une courbe $C$ associe la classe d'isomorphisme stable du
faisceau $\cN$ (resp. $\cE$) et Rao a montr\'e que cette application est surjective
et que ses fibres sont donn\'ees par le
th\'eor\`eme suivant (cf. aussi [N]) :

\th {Th\'eor\`eme de Rao, deuxi\`eme forme}. Soient $C,C'$ deux courbes et $\cN, \cN'$
(resp.
$\cE,\cE'$) les faisceaux intervenant dans une r\'esolution de type \sN\ (resp.
\sE) de $C$ et $C'$. Alors $C$ et $C'$ sont dans la m\^eme classe de biliaison si et
seulement si
$\cN$ et $\cN'$ (resp. $\cE$ et $\cE'$) sont stablement isomorphes, \`a d\'ecalage
pr\`es.

Il est facile de passer d'une forme \`a l'autre des th\'eor\`emes de Rao. Dans un sens
on utilise l'\'egalit\'e
$M_C= H^1_* \cN$ et dans l'autre le fait que
$\cN$ se retrouve comme faisceau associ\'e \`a un module de syzygies de la
r\'esolution libre de
$M_C$ (cf. [MDP1] II). En fait, on voit ais\'ement que le foncteur qui \`a un faisceau
$\cN$ (v\'erifiant
$H^2_* \cN =0$), \`a isomorphisme stable pr\`es, associe $H^1_* \cN$, \`a
isomorphisme pr\`es, est une \'equivalence de cat\'egories (cf. [Ho] 9.1). 

\vskip 0.3 cm

Nous allons montrer ici l'analogue du th\'eor\`eme de Rao, deuxi\`eme forme, pour les
familles (plates) de courbes de $\bP^3$ (param\'etr\'ees par un anneau noeth\'erien
local
$A$).\footnote {$(^1)$}{Dans un travail ult\'erieur [HMDP3] nous donnerons la variante
qui g\'en\'eralise le th\'eor\`eme de Rao premi\`ere forme et qui n\'ecessite
l'introduction d'objets nouveaux : les triades. Dans un autre article [HMDP2] nous
construirons des familles de courbes associ\'ees aux faisceaux localement libres,
c'est-\`a-dire l'analogue de la surjectivit\'e de l'application $C \mapsto \cN$.}  Pour
cela nous g\'en\'eralisons les notions de r\'esolutions de type
\sN\ et
\sE\ d'une famille
$\cC$, ce sont des suites exactes comme les suites (1) ci-dessus mais avec comme
seule hypoth\`ese le fait que $\cN$ et $\cE$ sont localement libres sur
$\bP^3_A$, nous montrons qu'il existe de telles r\'esolutions et que les faisceaux
$\cN$ et $\cE$ sont bien d\'etermin\'es  \`a  {\bf pseudo-isomorphisme}  pr\`es. Cette
notion de pseudo-isomorphisme s'exprime en termes de foncteurs de la cat\'egorie des
$A$-modules dans celle des $R_A$-modules gradu\'es : 

\th {D\'efinition 2.1}. Soient $\cN$ et $\cN'$  des faisceaux coh\'erents sur $\bP^3_A$ et
plats sur $A$  et soit
$f$ un morphisme de 
$\cN$ dans $\cN'$. On dit que $f$ est un pseudo-isomorphisme (en abr\'eg\'e un
\ps)  s'il induit~:\lign
0) un isomorphisme de foncteurs $H^0(\cN(n)\T_A \tas) \fl H^0(\cN'(n)\T_A \tas)$ pour
tout $n \ll 0$,  \lign
1) un isomorphisme de foncteurs $H^1_*(\cN\T_A \tas) \fl H^1_*(\cN'\T_A \tas)$ et
\lign 2) un monomorphisme de foncteurs $H^2_*(\cN\T_A \tas) \fl H^2_*(\cN'\T_A
\tas)$.\lign Deux faisceaux coh\'erents seront dits pseudo-isomorphes s'il existe
une cha\^\i ne de \ps\  qui les joint.

Cette notion g\'en\'eralise celle d'isomorphisme stable : on
montre en effet  (cf. 2.4) que
$f :
\cN
\fl
\cN'$ est un
\ps\ si et seulement si
$\cN$ est stablement isomorphe \`a une extension (pas n\'ecessairement scind\'ee) de
$\cN'$ par un dissoci\'e.  On obtient alors les th\'eor\`emes suivants, lorsque le
corps r\'esiduel de $A$ est infini :

\th {Th\'eor\`eme 3.1}.  Soient $\cC$ et $\cC'$ deux familles plates de courbes
param\'etr\'ees par l'anneau local  $A$. Alors, $\cC$ et
$\cC'$ sont dans la m\^eme classe de biliaison si et seulement si $\cJ_\cC$ et
$\cJ_{\cC'}$  sont pseudo-isomorphes, \`a d\'ecalage pr\`es.

\th {Th\'eor\`eme 3.2}. Soient $\cC$ et $\cC'$ deux familles plates de courbes
param\'etr\'ees par l'anneau local $A$, munies de r\'esolutions de type \sN\
(resp. \sE), avec des faisceaux $\cN$, $\cN'$ (resp. $\cE$, $\cE'$). Alors, $\cC$ et
$\cC'$ sont dans la m\^eme classe de biliaison si et seulement si $\cN$ et $\cN'$
sont pseudo-isomorphes,  \`a d\'ecalage pr\`es (resp. si $\cE^\vee$ et $\cE'^\vee$
sont pseudo-isomorphes,  \`a d\'ecalage pr\`es).

Pour avoir une variante de ces th\'eor\`emes qui ressemble plus \`a celle de Rao
on proc\`ede de la mani\`ere suivante : on note d'abord que toute classe de faisceaux
coh\'erents pour la relation de  pseudo-isomorphisme contient un faisceau
$\cN$ dit extraverti  c'est-\`a-dire v\'erifiant $\Ext^1_{R_A}
(H^0_*\cN, R_A)=0$. \footnote {$(^2)$}{Cette condition g\'en\'eralise la condition
$H^2_* (\cN) =0$, cf. 2.7.}  De plus ces faisceaux sont uniques dans leur classe, \`a
isomorphisme stable pr\`es. Toute famille de courbes admet une r\'esolution de type
\sN\ extravertie i.e. avec $\cN$ extraverti (resp. de type
\sE\ introvertie i.e. avec $\cE^\vee$ extraverti). On a alors :

\th {Th\'eor\`eme 3.3}. Soient $\cC$ et $\cC'$ deux familles plates de courbes
param\'etr\'ees par l'anneau local $A$, munies de r\'esolutions de type \sN\
extraverties (resp. de type \sE\ introverties), avec des faisceaux $\cN$, $\cN'$
(resp.
$\cE$,
$\cE'$). Alors,
$\cC$ et
$\cC'$ sont dans la m\^eme classe de biliaison si et seulement si $\cN$ et $\cN'$
(resp.  $\cE$ et $\cE'$)
sont stablement isomorphes \`a d\'ecalage pr\`es.

\vskip 0.3 cm

Pour une caract\'erisation de la liaison impaire,
cf. 3.10.
\vskip 0.3 cm

Les r\'esultats \'etablis ici dans le cas des courbes de $\bP^3$ se g\'en\'eralisent
sans difficult\'e au cas des sous-sch\'emas de codimension $2$ de $\bP^n$. On obtient
ainsi le th\'eor\`eme suivant (cf. 3.11 pour la d\'efinition de pseudo-isomorphisme) :

\th {Th\'eor\`eme 3.13}. Soient $\cC$ et $\cC'$ deux familles plates de sous-sch\'emas
de $\bP^n$ ($n \geq 3$),
de pure codimension $2$, sans composantes immerg\'ees, 
param\'etr\'ees par l'anneau local
$A$, munies de r\'esolutions de type
\sN\ (resp. \sE), avec des faisceaux $\cN$, $\cN'$ (resp. $\cE$, $\cE'$). Alors, $\cC$ et
$\cC'$ sont dans la m\^eme classe de biliaison si et seulement si $\cN$ et $\cN'$
sont pseudo-isomorphes,  \`a d\'ecalage pr\`es (resp. si $\cE^\vee$ et $\cE'^\vee$
sont pseudo-isomorphes,  \`a d\'ecalage pr\`es).

\titre {0. Rappels et notations}

Soit $A$  un anneau noeth\'erien. On pose $T = \Spec A$. On note $\bP^3_A$ ou
$\bP^3_T$ l'espace projectif de dimension $3$ sur $A$ et  $R_A$ l'anneau $A[X,Y,Z,T]$.
Si
$\cF$ est un faisceau coh\'erent sur
$\bP^3_A$ on note
 $H^i \cF$ le $A$-module $ H^i(\bP^3_A, \cF)$ et on pose $H^i_* \cF=\bigoplus _{n\in
\bZ} H^i \cF (n)$.
\vskip 0.3 cm

Une g\'en\'eralisation de la notion de module libre est celle de $R_A$-module {\bf 
dissoci\'e}. Il s'agit des $R_A$-modules gradu\'es  de la forme
$$F = \bigoplus_{i=1}^r M_i \T_A  R_A(-n_i), $$
o\`u les $ n_i$ sont des entiers et les $M_i$ des $A$-modules projectifs de type fini.
Un $R_A$-module dissoci\'e est projectif sur $R_A$ et m\^eme libre si $A$ est local.

De m\^eme, on introduit la notion de {\bf faisceau dissoci\'e}
sur $\bP^3_A$ : c'est un faisceau de la forme 
$$ \cF = \bigoplus_{i=1}^r M_i \T_A \cO_{\bP_A}(-n_i) $$
o\`u les $ n_i$ sont des entiers et les $M_i$ des $A$-modules projectifs de type fini. 
Si le module $F$ est dissoci\'e le faisceau associ\'e  $\widetilde F$ l'est aussi ; si
le faisceau
$\cF$ est dissoci\'e le module $H^0_* \cF$ l'est aussi.

Deux faisceaux $\cF$ et $\cF'$ sur $\bP^3_A$ sont dits {\bf stablement isomorphes}
s'il existe des faisceaux dissoci\'es $\cL$ et $\cL'$ et un isomorphisme $\cF \oplus
\cL
\simeq 
\cF' \oplus \cL'$. Cette relation est une relation d'\'equivalence.

\vskip 0.3 cm

Sur un corps $k$, on appelle courbe un sous-sch\'ema de $\bP^3_k$ de dimension $1$,
sans composante ponctuelle (immerg\'ee ou non), c'est-\`a-dire localement de
Cohen-Macaulay. On appelle surface un sous-sch\'ema \'equidimensionnel de dimension
$2$, donc un diviseur sur $\bP^3_k$. Une surface est d\'efinie par une unique
\'equation $f \in H^0 (\bP^3_k, \cO_{\bP} (d))$.

Soit  $T$ un sch\'ema. On rappelle qu'une famille de courbes gauches  $\cC$  sur
$T$ (on dira  simplement une courbe de $\bP^3_T$), est un sous-sch\'ema
ferm\'e de
$\bP^3_T$, plat sur
$T$, et dont les fibres sont des courbes au sens pr\'ec\'edent. On note $\cO_{\cC}$
le faisceau structural de $\cC$ et
$\cJ_{\cC}$ le faisceau d'id\'eaux qui d\'efinit $\cC$ dans $\bP^3_T$.

 Une surface de $\bP^3_T$ est un sous-sch\'ema ferm\'e $X$ de $\bP^3_T$, plat sur
$T$ tel que, pour tout $t \in T$, $X_t$ est une surface de $\bP^3_{k(t)}$. Si le sch\'ema
$T$ est connexe, le degr\'e de $X_t$ est  un entier $d$ ind\'ependant de $t$.

\vskip 0.3 cm

 Soit $\cM$ un faisceau coh\'erent sur $\bP^3_T$, plat sur $T$. Rappelons que la
dimension projective de $\cM$ (not\'ee $\dip \cM$) est le plus petit entier $n$ tel
qu'il existe une r\'esolution 
$$0 \fl \cP_n \fl \cP_{n-1} \fl \cdots \fl \cP_0 \fl \cM \fl 0$$
avec les $\cP_i$ localement libres. On a alors le r\'esultat suivant :

\th {Proposition 0.1}. Soit $A$ un anneau noeth\'erien, $T = \Spec A$ et soit 
 $\cM$ un faisceau coh\'erent sur $\bP^3_T$, plat sur $A$. Si $t$ est un point de $T$
on pose $\cM (t) = \cM\T_A k(t)$. Soit
$q$ un entier
$\geq 0$.  Les conditions suivantes sont
\'equivalentes~: \lign
i) $\dip \cM \leq q$, \lign
ii) $\cE xt^i_{\cO_\bP} (\cM, \cO_\bP) =0$ pour tout $i > q$, \lign
iii) pour tout point  $t \in T$, $\dip \cM (t) \leq q$, \lign
iv) pour tout point $t \in T$, ${\rm prof}\, \cM(t) \geq 3-q$. 

\dem L'\'equivalence de {\sl i)} et {\sl ii)} est standard, celle de {\sl i)} et de {\sl
iii)} r\'esulte de la platitude de $\cM$, enfin, celle de {\sl iii)} et {\sl iv)}
d\'ecoule du th\'eor\`eme d'Auslander-Buchsbaum.

\remas {0.2}  Notons deux cons\'equences de la proposition 0.1 : \lign
1) Le faisceau  $\cM$ est localement libre (i.e.  $\dip \cM = 0$) si et seulement si
on a 
$$\cE xt^i_{\cO_\bP} (\cM, \cO_\bP) =0 \quad \hbox {pour tout } \; i > 0.$$ 
2) En vertu de {\sl iv)},  la condition $\dip \cM \leq 2$ est \'equivalente \`a la
nullit\'e du foncteur 
$H^0
(\cM (n) \T_A \tas)$  pour $n \ll 0$. En effet, par dualit\'e (cf. [H] 7.4), cette
derni\`ere condition est
\'equivalente \`a $\Ext^3_{\cO_\bP} ( \cM , \omega_{\bP}(n) \T_A \tas) =\Ext^3_{\cO_\bP} (
\cM , \cO_{\bP}(n-4) \T_A \tas) =0$ pour $n \gg 0$, c'est-\`a-dire \`a la nullit\'e du
faisceau $\cE xt^3_{\cO_\bP} ( \cM, \cO_\bP)$. \lign 3) Si $\cC$ est une famille de
courbes param\'etr\'ee par $A$ on a $\dip \cO_\cC = 2$ et $\dip \cJ_\cC = 1$.

\titre {1. Liaison sur une base}

Soit $A$ un anneau noeth\'erien, $T = \Spec A$. On suppose $T$ connexe. \footnote
{$(^3)$}{Les r\'esultats de ce paragraphe sont valables sur un sch\'ema $T$
localement noeth\'erien.}

\th {D\'efinition 1.1}. Une courbe $\cC$ de $\bP^3_T$ est une intersection
compl\`ete de bidegr\'e $(d_1,d_2)$ si, pour tout $t \in T$, $\cC_t$ est une
intersection compl\`ete de bidegr\'e $(d_1,d_2)$ de $\bP^3_k$.

Attention,  sur un anneau quelconque la courbe $\cC$ est seulement localement sur
$T$  intersection de deux surfaces de $\bP^3_T$. 
Sur un anneau local on a la proposition suivante :

\th {Proposition 1.2}. On suppose $T = \Spec A$ o\`u $A$ est un anneau local
noeth\'erien. Soit $\cC$ une intersection compl\`ete de bidegr\'e $(d_1,d_2)$ de
$\bP^3_T$. On a une suite exacte :
$$ 0 \fl \cO_{\bP_A^3} (-d_1-d_2) \fl \cO_{\bP_A^3} (-d_1) \oplus \cO_{\bP_A^3}
(-d_2) \fl
\cJ_\cC \fl 0.$$

\th {D\'efinition 1.3}. Soit $\cD$ une intersection compl\`ete de $\bP^3_T$ et
soient 
$\cC_1$ et
$\cC_2$ deux courbes de $\bP^3_T$ sch\'ematiquement contenues dans $\cD$. Soit $p_i:
\cO_\cD
\fl
\cO_{\cC_i}$ la surjection canonique, qui induit un homomorphisme injectif 
$\f_i : \cH om_{\cO_\bP} (\cO_{\cC_i}, \cO_\cD) \fl \cO_\cD$. \lign
On dit que les courbes $\cC_1$ et
$\cC_2$ sont  {\bf li\'ees}
par 
$\cD$ si, pour $i=1,2$ (modulo $2$), l'image de $\f_i$ est l'id\'eal
$\cJ_{\cC_{i+1},\cD}$. En particulier on a les isomorphismes :
$$\cJ_{\cC_1,\cD} \simeq \cH om_{\cO_\bP} (\cO_{\cC_2}, \cO_\cD) \quad {\rm et} \quad
\cJ_{\cC_2,\cD} \simeq
\cH om_{\cO_\bP} (\cO_{\cC_1}, \cO_\cD).$$
On dit alors que l'on  passe de $\cC_1$ \`a $\cC_2$ par une liaison
\'el\'ementaire.

\vskip 0.3 cm

La proposition suivante est due \`a Kleppe, cf. [K] 2.4.

 \th {Proposition 1.4}. Soit $\cC$ une courbe de $\bP^3_T$ contenue dans une
intersection compl\`ete $\cD$. Soit $\f : \cH om_{\cO_\bP} (\cO_{\cC}, \cO_\cD) \fl
\cO_\cD$ l'homomorphisme canonique associ\'e. Alors l'image de $\f$ est l'id\'eal
$\cJ_{\cC',\cD} $ d'une courbe
$\cC'$ et les courbes $\cC$ et $\cC'$ sont li\'ees par
$\cD$. \lign
De plus, la liaison commute aux fibres, i.e., pour tout $t \in T$, les courbes
$\cC(t)$ et $\cC'(t)$ sont li\'ees par $\cD(t)$.

\th {D\'efinition 1.5}. Soient $\cC_1$ et $\cC_2$ deux courbes de $\bP^3_T$. On dit
qu'on passe de $\cC_1$ \`a $\cC_2$ par une biliaison \'el\'ementaire $(d,h)$ s'il
existe une surface $\cQ$ de $\bP^3_T$ de degr\'e $d$ et un entier $h$ tels que l'on
ait
$\cJ_{\cC_1,
\cQ}
\simeq
\cJ_{\cC_2, \cQ}(h)$.

La proposition suivante montre comment construire des biliaisons \'el\'ementaires :

\th {Proposition 1.6}. Soit $\cC$  une  courbe de $\bP^3_T$ et soient $d,h$ deux
entiers, avec
$d >0$.  Soit $\cQ$ une  surface (plate sur $A$) contenant $\cC$  et soit 
$u : \cJ_{\cC/\cQ}(-h) \fl \cO_{\cQ}$ un homomorphisme  tel que, dans
chaque fibre en $t \in T$, $u(t)$ soit injectif et non surjectif. Soit $\cO_{\cC'}$ le
conoyau de
$u$. Alors $\cC'$ est une famille (plate) de courbes et $u$ induit un isomorphisme de
$
\cJ_{\cC/\cQ}(-h)$ sur $
\cJ_{\cC/\cQ'}$,  de sorte que $\cC'$ est  obtenue
\`a partir de 
$\cC$ par une  biliaison
\'el\'ementaire  $(d,h)$. \lign Si $A$ est local il suffit de v\'erifier la
condition ci-dessus au point ferm\'e. \lign
Si
$h$ est
$\geq 0$ et si
$A$ est local on obtient de telles  biliaisons
\'el\'ementaires  (dites tri\-viales) en prenant pour $u$ la restriction \`a $
\cJ_{\cC/\cQ} $ de la multiplication par un \'el\'ement $H\in H^0
(\bP^3_T,
\cO_{\bP} (h))$, non diviseur de $0$ dans la fibre ferm\'ee de $\cQ$, c'est \`a
dire corres\-pondant
\`a une surface qui coupe proprement
$\cQ$ dans cette fibre. 

\dem Cela r\'esulte de [MDP1] III 2.3 et 2.6 et de [EGA] IV, 11.

\th {D\'efinition 1.7}.
On dit que $\cC_1$ et $\cC_2$ sont li\'ees (resp. bili\'ees) si, pour tout $t \in T$,
il existe un voisinage ouvert $U$ de $t$ tel que $\cC_1$ et $\cC_2$ sont jointes par
une suite de liaisons (resp. biliaisons) \'el\'ementaires. La relation de liaison (resp.
de biliaison) est une relation d'\'equivalence.

\th {Proposition 1.8}. Deux courbes $\cC_1$ et $\cC_2$ sont bili\'ees si et
seulement si elles sont jointes, localement sur $T$, par une suite form\'ee d'un nombre
pair de liaisons
\'el\'ementaires.

\dem 
On montre, exactement comme dans [MDP1] III 2.3, qu'une biliaison \'el\'ementaire sur
la surface $\cQ$  est
\'equivalente \`a deux liaisons \'el\'ementaires sur les surfaces $\cQ, \cS$ et $\cQ,
\cS'$. Il reste \`a voir que si, partant d'une courbe
$\cC$,  on effectue une liaison sur des surfaces $\cQ, \cS$ pour obtenir $\cC'$ puis 
une autre sur les surfaces
$\cQ',
\cS'$ pour obtenir $\cC''$ on peut passer de $\cC$ \`a $\cC''$ par des biliaisons
\'el\'ementaires. Pour cela on choisit une surface $\cT$, contenant $\cC'$, de degr\'e
assez grand pour que $\cT$ coupe proprement $\cQ$ et $\cQ'$ et on effectue
successivement les liaisons suivantes :
$$\cC \ \Fl {\cQ, \cS}\  \cC' \ \Fl {\cQ, \cT}\  \cC_0  \ \Fl {\cQ, \cT}\  \cC' \ \Fl
{\cQ',
\cT}\ 
\cC'_0\  \Fl {\cQ', \cT}\  \cC'\  \Fl {\cQ', \cS'}\  \cC''$$
et on passe ainsi de $\cC$ \`a $\cC''$ par trois biliaisons \'el\'ementaires sur les
surfaces $\cQ$, $\cT$, $\cQ'$.

\titre {2. Pseudo-isomorphismes et r\'esolutions de type \sE\ et \sN}

\tarte {a) D\'efinition des pseudo-isomorphismes}

\th {D\'efinition 2.1}. Soient $\cN$ et $\cN'$ des faisceaux coh\'erents sur $\bP^3_A$ et
plats sur $A$  et soit 
$f$ un morphisme de 
$\cN$ dans $\cN'$. On dit que $f$ est un {\bf  pseudo-isomorphisme} (en abr\'eg\'e un
\ps) s'il induit  :\lign
0) un isomorphisme de foncteurs $H^0(\cN(n)\T_A \tas) \fl H^0(\cN'(n)\T_A \tas)$ pour
tout $n \ll 0$,
\lign  1) un isomorphisme de foncteurs $H^1_*(\cN\T_A \tas) \fl H^1_*(\cN'\T_A \tas)$ et
\lign 2) un monomorphisme de foncteurs $H^2_*(\cN\T_A \tas) \fl H^2_*(\cN'\T_A
\tas)$. \lign  Deux faisceaux coh\'erents sur $\bP^3_A$ et plats sur $A$ seront dits
pseudo-isomorphes s'il existe une cha\^\i ne de \ps\  qui
les joint :
$$ \cN= \cN_0 \fl \cN_1 \lf \cN_2 \fl \cN_3 \lf \cdots \fl \cN_{2p-1} \lf \cN_{2p} =
\cN'.$$

\remas {2.2} \lign
0) La condition 0) de 2.1 est automatiquement v\'erifi\'ee si les faisceaux $\cN$ et
$\cN'$ sont localement libres, ou si ce sont des faisceaux d'id\'eaux ou si, plus
g\'en\'eralement, on a $\dip \cN \leq 2$ et $\dip \cN' \leq 2$. Dans tous
ces cas en effet
$H^0(\cN(n)\T_A
\tas)$ et
$H^0(\cN'(n)\T_A
\tas)$ sont nuls pour tout $n \ll 0$, cf. 0.2. Dans la suite de cet article nous serons la
plupart du temps dans cette situation.\lign 1) On peut inverser le sens des fl\`eches
dans la cha\^\i ne qui joint
$\cN$ et
$\cN'$. \lign 2)  Le compos\'e de
deux \ps\  en est un autre.  \lign
3) Si
$\cN$ est un faisceau coh\'erent et si
$\cL$ est dissoci\'e, l'injection canonique $\cN \fl \cN \oplus \cL$ et la projection
canonique
$\cN
\oplus \cL \fl \cN$ sont des \ps. Deux
faisceaux stablement isomorphes sont pseudo-isomorphes. \lign
4) La relation de pseudo-isomorphisme est une relation d'\'equivalence. 
\vskip 0.3 cm

Pour des faisceaux localement libres, on peut caract\'eriser les  \ps\  par
dualit\'e :

\th {Proposition 2.3}. Soient $\cN$ et $\cN'$  des faisceaux localement libres sur
$\bP^3_A$  et soit 
$f$ un morphisme de 
$\cN$ dans $\cN'$. Alors, $f$ est un \ps\  si et seulement si $f^\vee :
\cN'^\vee \fl \cN^\vee$  induit
\lign 1')  un isomorphisme de foncteurs $H^2_*(\cN'^\vee\T_A \tas) \fl H^2_*(\cN^\vee\T_A
\tas)$ et \lign 2') un \'epimorphisme de foncteurs $H^1_*(\cN'^\vee\T_A \tas) \fl
H^1_*(\cN^\vee\T_A
\tas)$.

\dem Cela r\'esulte de [H] 7.4 et de [AG] III 6.7.

\tarte {b) Caract\'erisation des \ps}

La caract\'erisation suivante des \ps\  est fondamentale  :

\th {Proposition 2.4}. Soit $f : \cN \fl \cN'$  un morphisme
de faisceaux coh\'erents sur $\bP^3_A$ et plats sur $A$. Alors,
$f$ est un
\ps\
 si et seulement si il existe  un faisceau dissoci\'e
$\cL$ et un morphisme $p : \cL \fl \cN'$ tel que l'on ait une suite exacte 
$$0 \fl \cS \fl \cN \oplus \cL \,\Fl {(f,p)} \, \cN' \fl 0$$
avec $\cS$ dissoci\'e.

\dem Si on a  suite exacte comme ci-dessus, il est clair que $f$ est un \ps\
(cela r\'esulte du fait que si $\cL$ est dissoci\'e on a $H^i_*(\cL \T \tas)=0$ pour
$i=1,2$ et $H^0_*(\cL(n) \T \tas)=0$ pour $n \ll 0$).

Inversement, supposons que $f : \cN \fl \cN'$ est un \ps. Il existe un entier $n_0$ tel
que, pour $n < n_0$, le morphisme $H^0\cN(n) \fl H^0\cN'(n)$
est un isomorphisme. Comme
$\cN'$ est coh\'erent, le module $N' = \bigoplus_{n \geq n_0} H^0\cN'(n)$ est de type
fini sur $R_A$. Soit $L$ un $R_A$-module libre de type fini avec une surjection
$L \fl N'$ qui donne, sur les faisceaux associ\'es  un morphisme surjectif
$p :
\cL
\fl
\cN'$.  Soit $\cS$ le noyau de $(f,p)
:
\cN
\oplus \cL \fl \cN'$. Il suffit de montrer que $\cS$ est dissoci\'e. 

Montrons d'abord que le morphisme de foncteurs 
$H^0_* ((\cN \oplus \cL) \T \tas) \fl H^0_* ( \cN' \T \tas)$ est un \'epimorphisme. Soit
$Q$ un $A$-module que l'on r\'esout par $0 \fl G \fl F \fl Q \fl 0$ avec $F$ libre sur
$A$. On a le diagramme commutatif de suites exactes :
$$ \matrix {H^0_* ((\cN \oplus \cL) \T F)\!\!&\!\!\fl\!\! & \!\!\!H^0_* ((\cN \oplus
\cL)
\T Q)\!&\!\fl\!& \!\!H^1_* ((\cN
\oplus \cL) \T G)\!&\!\fl\! &\!\!H^1_* ((\cN \oplus \cL) \T F)\cr
\Vf{a}&&\Vf{b}&& \Vf{c}&& \Vf{d} \cr
H^0_*(\cN' \T F)& \fl&H^0_*(\cN' \T Q )&\fl&H^1_*(\cN' \T G)&\fl&H^1_*(\cN' \T F)\cr }$$
Comme $f$ est un \ps, les morphismes $c$ et $d$ sont des isomorphismes. Par ailleurs,
comme $F$ est libre, 
$a$ est surjectif (pour $n < n_0$ parce que $H^0\cN(n) \fl H^0\cN'(n)$ l'est, pour $n \geq
n_0$ parce que
$H^0\cL(n) \fl H^0\cN'(n)$ l'est). Il en r\'esulte que
$b$ est surjectif, comme annonc\'e.

Alors, le morphisme de foncteurs 
$H^i_* ((\cL\oplus\cN) \T \tas) \fl H^i_* (\cN' \T \tas) $ est un \'epimorphisme pour
$i=0$, un isomorphisme pour $i=1$ et un monomorphisme pour $i=2$   et il en r\'esulte que
$H^i_* (\cS \T \tas) $ est nul pour $i= 1,2$. De plus, comme on a
$H^0 (\cL(n)\T \tas)=0$ pour $n\ll 0$, il r\'esulte de la condition 0) de la d\'efinition
des
\ps\ que l'on a $H^0 (\cS(n)\T \tas)=0$ pour $n \ll0$. Cela montre que $\cS$ est
dissoci\'e en vertu de [H] 7.9.

\remas {2.5} \lign
1) Une d\'emonstration analogue \`a celle effectu\'ee ci-dessus montre que dans la
d\'efinition 2.1 on peut remplacer la condition 0) par le fait que, pour $n \ll 0$, le
morphisme
$H^0\cN(n)
\fl H^0\cN'(n)$ est un isomorphisme.
\lign
2)  Soit $f : \cN \fl \cN'$ un \ps\  de faisceaux
coh\'erents.  La ca\-rac\-t\'erisation 2.4 montre que si $\cN'$ est localement libre il
en est de m\^eme de
$\cN$.  En revanche, $\cN$ peut \^etre localement libre sans que $\cN'$ le soit, cf.
2.18.1. Cependant on a $\dip \cN \leq 1 \equi \dip \cN' \leq 1$ car ces conditions sont
\'equivalentes \`a la nullit\'e de
$\cExt^i_{\cO_\bP}
(\cN,
\cO_\bP)$ pour $i \geq 2$, cf. 0.1.

\tarte {c) Faisceaux extravertis}

\th {Proposition-D\'efinition 2.6}. Soit $\cN$ un faisceau coh\'erent sur
$\bP^3_A$ et plat sur $A$. Posons $N=H^0_*\cN$.  Les conditions suivantes sont
\'equivalentes : \lign
 i) $\Ext^1_{R_A} (N, R_A)=0$, \lign
ii) $\Ext^1_{\cO_\bP} (\cN, \cO_\bP(*))=\bigoplus_{n \in \bZ}\Ext^1_{\cO_\bP} (\cN,
\cO_\bP(n))=0$,
\lign iii) $ H^1_* \cN^\vee = 0$ et $\cE xt^1_{\cO_\bP} (\cN, \cO_\bP) =0$. \lign
Si on suppose, de plus, qu'on a $\dip \cN \leq 1$, ces conditions sont aussi
\'equivalentes \`a : \lign
iv) $ H^1_* \cN^\vee = 0$ et $\cN$ localement libre. \lign
 Un faisceau
$\cN$  coh\'erent sur
$\bP^3_A$ et plat sur $A$ est dit {\bf extraverti}  s'il v\'erifie les conditions i) ii)
iii) ci-dessus. 

\dem Pour montrer l'\'equivalence de {\sl i} et {\sl ii} il suffit d'\'etablir
 l'\'egalit\'e :
$$\Ext^i_{R_A} (N, R_A)=\Ext^i_{\cO_\bP} (\cN, \cO_\bP(*)) \qquad  {\rm pour }\
i=0,1,2.$$
Cette relation est claire pour $i=0$ et on obtient les autres cas par r\'ecurrence en
uti\-li\-sant une suite exacte $0 \fl M \fl P \fl N \fl 0$ avec $P$ libre et en tenant
compte de l'\'egalit\'e
$\Ext^i_{\cO_\bP} (\cP, \cO_\bP(*))=\Ext^i_{R_A} (P, R_A)=0$ pour $i=1,2$.

Pour  l'\'equivalence de {\sl ii} et {\sl iii} on utilise la suite spectrale des Ext
dont la suite exacte des termes de bas degr\'e s'\'ecrit :
$$0 \fl H^1_* \cN^\vee \fl \Ext^1_{\cO_\bP} (\cN, \cO_\bP(*)) \fl H^0_* (\cE
xt^1_{\cO_\bP} (\cN, \cO_\bP)) \fl H^2_* \cN^\vee $$
Il est alors clair que {\sl iii} implique {\sl ii} et que {\sl ii} entra\^\i ne la
nullit\'e de $H^1_* \cN^\vee$. De plus, pour $n \gg 0$, on a $H^i \cN^\vee (n) =0$ pour
$i=1,2$, donc aussi $ H^0 (\cE
xt^1_{\cO_\bP} (\cN, \cO_\bP(n))) =0$ ce qui implique la nullit\'e du faisceau $\cE xt$. 

Le cas $\dip \cN \leq 1$ r\'esulte de 0.2.1.

\rema {2.7} La condition $\cN$ extraverti g\'en\'eralise la condition usuelle
$H^2_* \cN =0$ impos\'ee aux faisceaux (localement libres) des r\'esolutions de type \sN\
sur un corps. En effet, si on pose
$C= H^2_* (\cN)$ on a, par dualit\'e, $ H^1_* \cN^\vee = C^*(4)$ et
$C$ est nul si et seulement si $C^*$ l'est, de sorte que la condition donn\'ee est
bien \'equivalente \`a  la nullit\'e de $H^2_* \cN$.

\vskip 0.3 cm

Dans le cas des faisceaux extravertis les \ps\ ne sont rien d'autres que
les isomorphismes stables :

\th {Corollaire 2.8}.  Soit $f : \cN \fl \cN'$ un morphisme
de faisceaux  coh\'erents. On suppose $\cN'$ extraverti.
Alors,
$f$ est un
\ps\ si et seulement si $\cN$ et $\cN'$ sont stablement
isomorphes {\it via} $f$.

\dem Supposons que $f$ est un \ps\ (l'autre sens est \'evident avec 2.2.3). En vertu de
 2.4
on a une suite exacte
$0\fl\cS\fl\cN\oplus\cL\,\fl\,\cN'\fl 0$ et cette
suite est scind\'ee parce que $\Ext^1_{\cO_\bP} (\cN', \cS)=0$ (c'est vrai
pour $\cS = \cO_\bP(n)$ puisque $\cN'$ est extraverti, donc aussi pour
$\cS$ dissoci\'e).

\vskip 0.3 cm

Le lemme suivant va permettre de construire des faisceaux extravertis :

\th {Lemme 2.9}. Soit $M$ un module de type fini sur $R_A$. Il
existe un  module
$P$ dissoci\'e, un module de type fini $N$ qui v\'erifie $\Ext^1_{R_A}(N,R_A)=0 $
 et une suite exacte : $0\to P\to N \to M \to 0$. 

\dem  On consid\`ere une couverture minimale $P'\to
\Ext^1_{R_A}(M,R_A)$, o\`u $P'$ est un $R_A$-module libre gradu\'e. La
composition des fl\`eches :
$$R_A\to P'\otimes P^{'\vee}\to \Ext^1_{R_A}(M,R_A)\otimes P^{'\vee}\to
\Ext^1_{R_A}(M, P^{'\vee})$$
fournit une extension de $M$ par  $P^{'\vee}$, donc une suite exacte $0\to
P^{'\vee}\to
N \to M \to 0$. Lorsqu'on la dualise, on obtient un cobord $P'\to
\Ext^1_{R_A}(M,R_A)$ qui est l'homomorphisme dont on est parti, donc qui est
surjectif. Autrement dit, on a $\Ext^1_{R_A}(N,R_A)=0$, et en posant
$P=P^{'\vee}$, on a la suite annonc\'ee.

\th {Proposition 2.10}. \lign 
1) Soit $\cM$ un faisceau coh\'erent sur $\bP^3_A$ et plat sur $A$. Il existe un
faisceau  extraverti $\cN$  et un faisceau dissoci\'e $\cP$ et une suite exacte $0 \fl
\cP \fl \cN \Fl {p} \cM \fl 0$, de sorte que $p$ est un \ps.
\lign 2) Si on suppose $\dip \cM \leq 1$ (par exemple si $\cM$ est localement libre sur
$\bP^3_A$ ou si $\cM= \cJ_\cC$ est le faisceau d'id\'eaux d'une famille de courbes) le
faisceau 
$\cN$ est localement libre.

\dem On consid\`ere  le $R_A$-module de
type fini $M = \bigoplus _{n \geq 0} H^0 \cM (n)$. Le point 1)
r\'esulte aussit\^ot de 2.9 appliqu\'e \`a $M$. 

Pour le point 2) on note  que l'hypoth\`ese $\dip \cM \leq 1$
impose $\cExt^i_{\cO_\bP} (\cM,
\cO_\bP) =0$ pour $i \geq 2$, donc encore
$\cExt^i_{\cO_\bP} (\cN,
\cO_\bP) =0$ pour $i \geq 2$. Comme $\cN$ est extraverti on a aussi
$\cExt^1_{\cO_\bP} (\cN,
\cO_\bP)=0$  en vertu de 2.6 et il en r\'esulte que  $\cN$ est localement
libre (cf. 0.2.1).

\tarte {d) Le lemme de Verdier}

Nous montrons ici un lemme qui permet de r\'eduire la longueur des
cha\^\i nes de \ps\ \`a deux. C'est un analogue du calcul de fractions introduit
par Verdier dans [V].

\th {Proposition 2.11}. (Lemme de Verdier) Soient $\cN'$ et $\cN''$ 
des faisceaux coh\'erents et plats sur $A$ (resp.
localement libres sur $\bP^3_A$) pseudo-isomorphes. Il existe un faisceau coh\'erent et
plat sur $A$ (resp. localement libre sur $\bP^3_A$) $\cN$ et des
\ps\  : 
$\cN' \lf \cN \fl \cN''$. De plus, on peut supposer $\cN$ extraverti.

\dem  Le r\'esultat
provient, par r\'ecurrence sur la longueur d'une cha\^\i ne de \ps\ joignant les
faisceaux, du lemme suivant~:

\th {Lemme 2.12}. Soient $\cN, \cN', \cN''$ des faisceaux coh\'erents et plats sur
$A$ (resp. localement
libres sur $\bP^3_A$) avec des
\ps\
$f': \cN' \fl \cN$ et $f'': \cN'' \fl \cN$. Il existe un faisceau coh\'erent plat
sur $A$ (resp. localement
libre sur $\bP^3_A$)
$\cM$ et des \ps\  $g' : \cM \fl \cN'$ et $g'' : \cM \fl \cN''$.

\dem (du lemme 2.12) Vu la remarque 2.5.2, il suffit de traiter le cas coh\'erent.
Quitte
\`a ajouter
\`a 
$\cN'$ et
$\cN''$  des faisceaux dissoci\'es on peut supposer qu'on a des suites exactes
$0 \fl \cS' \fl \cN' \fl \cN \fl 0$ avec $\cS'$
dissoci\'e et l'analogue avec $\cN''$. On prend alors pour $\cM$ le produit fibr\'e
$\cM = \cN' \times_\cN \cN''$ et on a le diagramme de suites exactes suivant :
$$ \matrix { &&0&&0&&0\cr
&&\vf&&\vf&&\vf\cr
0&\fl& \cS'&\fl&\cN'&\fl&\cN&\fl&0\cr
&&\vf&&\vf&&\parallel\cr
0&\fl& \cM&\fl&\cN'\oplus \cN''&\fl&\cN&\fl&0\cr
&&\Vf{g''}&&\vf&&\vf\cr
0&\fl& \cN''&=&\cN''&\fl&0\cr
&&\vf&&\vf\cr
 &&0&&0\cr
}$$ 
On a ainsi le \ps\ $g''$ cherch\'e et on proc\`ede de m\^eme pour $g'$.

L'assertion sur extraverti r\'esulte de 2.10.

\th {Proposition 2.13}. Soient $\cN$ et $\cN'$  des faisceaux
extravertis. Alors $\cN$ et $\cN'$ sont  pseudo-isomorphes  si et seulement si ils sont
stablement isomorphes.

\dem C'est clair avec 2.11, 2.8 et le fait que la relation d'isomorphisme stable est une
relation d'\'equivalence.

\vskip 0.3 cm

Dans le cas d'un anneau local, la structure de la classe d'isomorphisme stable est
particuli\`erement simple :

\th {Proposition 2.14}. On suppose $A$ local. La classe de pseudo-isomorphisme  d'un
faisceau $\cN$ coh\'erent et v\'erifiant $\dip \cN \leq 1$ contient un  faisceau  $\cN_0$
localement libre, extraverti, sans facteur direct dissoci\'e.  Tout autre faisceau
coh\'erent
 extraverti de la classe est de la forme $\cN_0 \oplus
\cL$ o\`u $\cL$ est dissoci\'e. Le faisceau $\cN_0$ est unique \`a isomorphisme pr\`es.

\dem L'existence est \'evidente avec 2.10.  Si $\cN$ est un autre faisceau extraverti
de la classe, il est stablement isomorphe \`a $\cN_0$ par 2.13. On montre alors qu'il
est de la forme $\cN_0 \oplus \cL$ en raisonnant comme dans [N] Proposition 2.3.
(Dans [N] l'assertion est donn\'ee dans le cas d'un corps, mais on la g\'en\'eralise
au cas d'un anneau local en notant que si on a un morphisme surjectif
$(f,f') : M
\oplus M'
\fl A$ de
$A$-modules, alors
$f$ ou
$f'$ est surjectif). L'unicit\'e en d\'ecoule aussit\^ot.

\tarte {e) R\'esolutions de type \sE\ et \sN}

\th {D\'efinition 2.15}. Soit $\cC$ une courbe de $\bP^3_T$. Une r\'esolution de type
\sN\ (resp. \sE) de $\cC$ est une  suite exacte de faisceaux sur $\bP^3_T$:
$$0 \fl \cP \fl \cN \fl \cJ_\cC \fl 0 \quad {\rm (resp.} \quad 0 \fl \cE \fl \cF\fl
\cJ_\cC
\fl 0)$$
o\`u $\cP$ et $\cF$ sont dissoci\'es  et o\`u $\cN$ et
$\cE$ sont  localement libres.

\th {D\'efinition 2.16}. Soient $r$ et $r'$ deux r\'esolutions de type \sN\ d'une
courbe $\cC$. 
Un morphisme $r \fl r'$ (resp. un isomorphisme) consiste en  la donn\'ee de deux
morphismes  de faisceaux
$u,v$ (resp. deux isomorphismes) rendant commutatif le diagramme suivant :
$$\matrix{ (r) &&&&&0 & \fl & \cP& \fl & \cN& \fl &\cJ_\cC& \fl & 0 \cr
\vf&&&&&&&\Vf{u}&&\Vf{v}&& \parallel \cr
(r')&&&&&0 & \fl & \cP'& \fl & \cN'& \fl &\cJ_\cC& \fl & 0 \cr
}$$
On a \'evidemment une d\'efinition analogue avec les r\'esolutions de type \sE.

\remas {2.17} \lign
0) On notera que  si on a une r\'esolution $0 \fl \cP \fl \cN \fl \cJ_\cC \fl 0$
(resp. $0 \fl \cE \fl \cF\fl
\cJ_\cC
\fl 0$), avec $\cN$ (resp. $\cE$) coh\'erent et $\cP$ (resp. $\cF$) dissoci\'e, alors
 $\cN$
  v\'erifie $\dip \cN \leq 1$ (resp. $\cE$  est localement libre). Cela
r\'esulte de 0.1  et de l'\'egalit\'e
$\cE xt^i_{\cO_\bP} (\cJ_\cC, \cO_\bP)=0$  pour $i \geq 2$.
\lign 1) Si
$0
\fl
\cP
\Fl {u}
\cN
\fl
\cJ_\cC
\fl 0
$ est une r\'esolution de type \sN\ et  $\cL$ un faisceau dissoci\'e quelconque on obtient
une autre r\'esolution de type \sN\ en prenant :
$$0 \fl \cP \oplus \cL \Fl {u \oplus 
{\rm Id}} \cN \oplus \cL \fl \cJ_\cC \fl 0.$$
On a, bien entendu, une assertion analogue pour les r\'esolutions de type \sE. \lign
2) Soit $r$  une r\'esolution de type \sN\ : $0 \fl \cP \fl \cN \fl \cJ_\cC
\fl 0$ et soit un $f: \cN' \fl \cN$
 un \ps. En vertu de 2.4, il existe des faisceaux dissoci\'es $\cS$ et $\cL$ et une suite
exacte 
$0 \fl \cS \fl \cN' \oplus \cL \Fl {(f,p)} \cN \fl 0$. On en d\'eduit une r\'esolution
$r'$ de type \sN\ :
$0 \fl \cP\oplus \cS \fl \cN'\oplus \cL \fl \cJ_\cC \fl 0$ et  un morphisme de
$r'$ dans
$r$ induit par $(f,p)$.
\lign 3) Si on
a des r\'esolutions de type
\sN\ et 
\sE\ de
$\cJ_\cC$ les faisceaux $\cN$ et $\cE$
sont reli\'es, avec les notations de 2.15, par la suite exacte suivante : $0 \fl \cE \fl
\cP
\oplus
\cF
\fl
\cN
\fl 0$ (passer aux modules et relever la fl\`eche de $F$ dans $I_\cC$ en une fl\`eche
de
$F$ dans $N$).

\tarte {f) Lien avec les \ps}

\th {Proposition 2.18}. \lign
1)  Soit $\cC$ une courbe de $\bP^3_A$ munie d'une
r\'esolution de type \sN\ : $0 \fl \cP \fl \cN \Fl{p} \cJ_\cC \fl 0$. Alors $p$ est
un \ps. \lign
2) Si $\cC$ est une courbe de $\bP^3_A$ et $\cQ$ une surface contenant $\cC$, le
morphisme $\cJ_\cC \fl \cJ_{\cC/\cQ}$ est un \ps. \lign
3)  Si on a deux r\'esolutions de type \sN\ de $\cC$ les
faisceaux $\cN$ intervenant dans ces r\'esolutions sont pseudo-isomorphes.
 Si ces
r\'esolutions sont extraverties  (cf. 2.22)
les faisceaux cor\-res\-pondants
sont stablement isomorphes.

\dem Dans les deux cas 1) et 2) cela r\'esulte de 2.4 : c'est clair pour 1) et, pour
2), si
$\cQ$ est de degr\'e $s$, on a la suite exacte $0 \fl \cO_\bP (-s) \fl  \cJ_\cC \fl
\cJ_{\cC/\cQ} \fl 0$.

Si on a deux r\'esolutions de type \sN, les faisceaux $\cN$ et $\cN'$ sont tous deux
pseudo-isomorphes \`a $\cJ_\cC$ en vertu de 2.4, ce qui donne la premi\`ere assertion de
3). Le cas extraverti  r\'esulte de 2.13.

\vskip 0.3 cm

Lorsqu'on a un morphisme de r\'esolutions le r\'esultat suivant pr\'ecise 2.18.3 et
vaut aussi pour les r\'esolutions de type \sE\ :

\th {Proposition 2.19}. Soit $\cC$ une courbe et $\f : r \fl r'$ un morphisme
entre deux r\'esolutions de type \sN\ (resp. \sE) de $\cC$. Alors le morphisme
induit par $\f$ sur les faisceaux $\cN \fl \cN'$ (resp. $\cE \fl \cE'$) est un \ps\
 (resp.
le dual d'un \ps).

\dem Faisons le pour le cas du type \sN, l'autre est analogue.  Cela r\'esulte du
fait que, comme $\cP$ est dissoci\'e, donc n'a ni $H^1$ ni $H^2$ on a des
isomorphismes des foncteurs $H^1_* (\cN \T \tas) \simeq H^1_* (\cJ_\cC \T \tas) \simeq
H^1_* (\cN' \T \tas)$ et un diagramme commutatif :
$$ \matrix {0& \fl & H^2_*(\cN \T \tas) & \fl & H^2_*(\cJ_\cC \T \tas)& \fl & \cdots \cr
&& \Vf{\f_2}&& \parallel \cr
0& \fl & H^2_*(\cN' \T \tas) & \fl & H^2_*(\cJ_\cC \T \tas)& \fl & \cdots \cr
}$$
qui montre que $\f_2$ est un monomorphisme de foncteurs.

\tarte {g) Existence des r\'esolutions de type \sE}

\th {Proposition 2.20}.  Soit $\cC$ une courbe de $\bP^3_A$. Il existe une
r\'esolution
$r$ de type \sE\ de $\cC$. On peut  supposer cette r\'esolution introvertie,
c'est-\`a-dire telle que $\cE^\vee$ soit extraverti.

\dem  On consid\`ere $I_\cC= H^0_* \cJ_\cC$. Soit
$F$ un
$R_A$-module libre avec
$p: F \fl I_\cC$ surjectif  et soit $E$ le noyau de $p$. Alors la suite de faisceaux
associ\'es $0 \fl \cE \fl \cF \fl \cJ_\cC \fl 0$ est une r\'esolution de type \sE\
introvertie.

En effet, comme $\cJ_\cC$ est plat sur $A$ il en est de m\^eme de $\cE$. De plus,
dans chaque fibre en $t \in T$, $\cE(t)$ est localement libre, donc $\cE$ est localement
libre. Par ailleurs, 
$\cE^\vee$ est extraverti car $p$ est surjective.

\th {Proposition 2.21}. (Domination)
Soient $r_1$ et $r_2$  deux r\'esolutions de type \sE\ de $\cC$. Il existe une
r\'esolution $r$ de type \sE\   introvertie avec des morphismes $r_1 \fl r$ et $r_2
\fl r$.  \lign
Les faisceaux  $\cE_1^\vee$ et $\cE_2^\vee$ correspondants  sont pseudo-isomorphes. Si ces
r\'esolutions sont  introverties les faisceaux $\cE_1$ et $\cE_2$  sont stablement
isomorphes.

\dem Posons, pour $k=1,2$, $(r_k) = (0 \fl \cE_k \fl \cF_k \fl \cJ_\cC \fl 0)$ et soit
$I_k$ l'image du module libre $F_k$ dans $H^0_* \cJ_\cC$ de sorte que  $I_1$ et $I_2$
sont contenus dans $I_\cC$. Si on prend un module libre $F$ avec une suite exacte
$0 \fl E \fl F \fl I_\cC \fl 0$ on v\'erifie aussit\^ot que les inclusions des
$I_k$ dans $I_\cC$ se rel\`event aux suites $ 0 \fl E_k \fl F_k \fl I_k \fl
0$ et fournissent les morphismes souhait\'es.

En vertu de 2.19 les faisceaux $\cE_1$ et $\cE_2$ sont pseudo-isomorphes. Le cas
introverti d\'ecoule de 2.13.

\tarte {h) Existence de r\'esolutions de type \sN}

\th {Proposition 2.22}.  Soit $\cC$ une courbe de $\bP^3_A$. Il existe une
r\'esolution de type \sN\ de $\cC$. On peut  supposer cette r\'esolution
extravertie c'est-\`a-dire telle que $\cN$ soit extraverti.

\dem Comme le faisceau $\cJ_\cC$ est plat sur $A$ et v\'erifie $\dip \cJ_\cC \leq 1$,
cela r\'esulte de 2.10.2.

\th {Proposition 2.23}. (domination) Soient $r_1$ et $r_2$  
deux r\'esolutions de type \sN\ de $\cC$. Il existe une
r\'esolution $r$  de type \sN\ extravertie  avec des morphismes $r\fl r_1$ et
$r\fl r_2$.

\dem En utilisant 2.11 et les remarques 2.17.1 et 2.17.2 on se ram\`ene au cas o\`u les
deux r\'esolutions sont extraverties avec le m\^eme $\cN$ et il s'agit de montrer qu'il y
a un morphisme de l'une dans l'autre. Comme on a $\Ext^1_{R_A} (N,R_A)=0$, donc
aussi $\Ext^1_{R_A} (N,P)=0$, puisque $P$ est dissoci\'e, l'identit\'e de $I_\cC =
H^0_*
\cJ_\cC$ se rel\`eve en une fl\`eche de $N$ dans $N$ qui donne le morphisme cherch\'e.

\tarte {i) Comportement des r\'esolutions par liaison}

Si $\cC$ est une courbe de $\bP^3_T$, le faisceau dualisant relatif de $\cC$ sur $T$
est, par d\'efinition :
$$\omega_{\cC/T} = \cExt ^2_{\cO_{\bP_T}} (\cO_{\cC}, \cO_{\bP_T}(-4))$$
et on a les r\'esultats suivants (cf. [S], 3.2 et 3.3, voir aussi [BPS] et [JS]) :

 $$\cExt ^i_{\cO_{\bP_T}} (\omega_{\cC/T}, \cO_{\bP_T}) = \cases {0& si $i\neq 2$,
\cr
\cO_{\cC}(4)& si $i=2$. \cr} \leqno {1)}$$
$$\cExt ^i_{\cO_{\bP_T}} (\cO_{\cC}, \cO_{\bP_T}) = \cases {0& si $i\neq 2$,
\cr
\omega_{\cC/T}(4)& si $i=2$. \cr} \leqno {2)}$$
\vskip 0.3 cm

Soit $\cC_1$ une courbe de $\bP^3_A$, contenue dans une intersection compl\`ete
$\cD$, de bidegr\'e $(s, t)$. On suppose qu'on a une r\'esolution de
$\cD$ comme en 1.2 (c'est le cas, par exemple, si $A$ est local) :
$$0 \fl \cO_\bP (-s-t) \fl \cO_\bP (-s) \oplus \cO_\bP (-t) \fl \cJ_\cD \fl 0.$$
Soit $\cC_2$ la courbe li\'ee \`a $\cC_1$ par $\cD$ (cf. 1.4). Les deux
propositions suivantes  comparent les r\'esolutions de $\cC_1$ et $\cC_2$ :

\th {Proposition 2.24}. Soit $0 \fl \cP \fl \cN \fl \cJ_{\cC_1} \fl 0 $ une
r\'esolution de type \sN\ de $\cC_1$. On a une r\'esolution de type \sE\ de $\cC_2$ :
$$0 \fl \cN^\vee (-s-t) \fl \cP^\vee (-s-t) \oplus \cO_{\bP} (-s) \oplus \cO_{\bP}
(-t) \fl \cJ_{\cC_2} \fl 0.$$

\dem La m\'ethode est classique, cf. [PS] 2.5. Comme la fl\`eche $H^0_* \cN \fl H^0_*
\cJ_{\cC_1}$ est surjective, l'inclusion
$\cJ_\cD \subset \cJ_{\cC_1}$ se rel\`eve en un morphisme des r\'esolutions. En
dualisant ce morphisme on obtient le morphisme de complexes suivant (en degr\'es $1$
et $0$) :
$$ \matrix { (\Lpont)\quad &&&\cN^\vee & \fl &\cP^\vee \cr
&&&\vf&& \vf\cr
(\Mpont) \quad && &\cO_\bP (s) \oplus \cO_\bP (t) & \fl& \cO_\bP (s+t) \cr}$$
Les groupes d'homologie en degr\'e $1$ sont tous deux isomorphes \`a $\cO_\bP$ et en
degr\'e $0$ on trouve respectivement $h_0 \Lpont =\cExt^1_{\cO_\bP}
(\cJ_{\cC_1},
\cO_\bP) \simeq \omega_{\cC_1/T} (4) $ et $h_0 \Mpont =\cExt^1_{\cO_\bP}
(\cJ_\cD,
\cO_\bP) = \omega_{\cD/T} (4)$ en vertu du rappel.

On consid\`ere la suite exacte $0 \fl \cJ_\cD \fl \cJ_{\cC_1} \fl \cJ_{\cC_1,\cD} \fl
0$. La d\'efinition de la liaison donne un isomorphisme : $\cJ_{\cC_1,\cD} \simeq \cH
om_{\cO_\bP} (\cO_{\cC_2}, \cO_\cD) $ et, vu la r\'esolution de $\cO_\cD$, ce faisceau
est encore isomorphe \`a  $\cExt^2_{\cO_\bP}
(\cO_{\cC_2},
\cO_\bP (-s-t)) = \omega_{\cC_2/T} (4-s-t)$. On a donc la suite exacte 
$$0 \fl \cJ_\cD \fl \cJ_{\cC_1} \fl \omega_{\cC_2/T} (4-s-t) \fl 0$$
qui, par dualit\'e donne 
$$0 \fl \omega_{\cC_1/T}Ê(4) \fl \omega_{\cD/T} (4) \fl \cO_{\cC_2} (s+t) \fl 0$$ en
vertu des r\'esultats de Schlesinger.

 On consid\`ere alors le c\^one $\Cpont$ du morphisme $\Lpont \fl \Mpont$ :
$$(\Cpont) \qquad \cN^\vee \fl \cP^\vee \oplus \cO_\bP (s) \oplus \cO_\bP (t)  \fl
\cO_\bP (s+t).$$ On a $h^i \Cpont =0$ pour $i \geq 1$ et la suite exacte d'homologie :
$ 0  \fl h_0 \Lpont \fl h_0 \Mpont \fl h_0 \Cpont \fl 0$
 montre qu'on a $h_0 \Cpont = \cO_{\cC_2} (s+t)$, d'o\`u la suite exacte $$0 \fl
\cN^\vee \fl \cP^\vee \oplus \cO_\bP (s) \oplus \cO_\bP (t)  \fl
\cO_\bP (s+t) \fl \cO_{\cC_2} (s+t) \fl 0$$
et la conclusion s'ensuit.

\th {Proposition 2.25}. Soit $0 \fl \cE \fl \cF \fl \cJ_{\cC_1} \fl 0 $ une
r\'esolution de type \sE\ de $\cC_1$. On suppose que le morphisme compos\'e $
\cO_{\bP} (-s)
\oplus \cO_{\bP} (-t) \fl \cJ_\cD \subset \cJ_{\cC_1}$ se rel\`eve \`a $\cF$. Alors
on a une r\'esolution de type
\sN\ de
$\cC_2$ :
$$0 \fl \cF^\vee (-s-t) \fl \cE^\vee (-s-t) \oplus \cO_{\bP} (-s) \oplus \cO_{\bP}
(-t) \fl \cJ_{\cC_2} \fl 0.$$

\dem L'existence du rel\`evement \'etant assur\'ee, la d\'emonstration est identique
\`a la pr\'ec\'edente.

\remas {2.26}  \lign 
1)  L'hypoth\`ese de l'existence d'un rel\`evement dans 2.25 est
v\'erifi\'ee si le faisceau $\cE$ v\'erifie $H^1_* \cE=0$, ou si on prend $s,t > n_0$ o\`u
$n_0$ est le plus petit entier au-del\`a duquel on a $H^1 \cE(n) =0$. \lign
2) On peut utiliser 2.25 pour donner des variantes de l'existence des r\'esolutions
de type \sN. \lign
3) Il r\'esulte de  2.25 et 2.26 et de 1.8 que, par biliaison, les faisceaux $\cE$ et
$\cN$ sont invariants, \`a adjonction d'un faisceau dissoci\'e pr\`es.

\titre {3. Le th\'eor\`eme de Rao}

\tarte {a) \'Enonc\'es}

Dans ce paragraphe on suppose l'anneau $A$  {\bf local} noeth\'erien et on note
$m_A$ son radical. 
Soient
$t$ le point ferm\'e de $T = \Spec A$ et $k(t)$ le corps r\'esiduel, qui sera suppos\'e
{\bf infini}.

Soient $H$  un $A$-module de type fini, $\ov H =H \T_A k(t)$ le $k(t)$-espace
vectoriel de
dimension finie obtenu par r\'eduction modulo $m_A$ et $\pi : H \fl \ov H$
la projection
canonique. Nous dirons qu'une propri\'et\'e $P$ des \'el\'ements de $H$ est
vraie pour $h$
``g\'en\'eral''dans $H$ s'il existe un ouvert de Zariski non vide
$\ov U $ du sch\'ema affine associ\'e \`a $\ov H$ tel que $P$ soit vraie pour tout 
$h \in \pi^{-1} (\ov U)$. Puisque $k(t)$ est infini,
 un tel ouvert a des points rationnels, et son image r\'eciproque $\pi^{-1} (\ov U)$,
consid\'er\'ee comme sous-ensemble du $A$-module $H$, n'est pas vide.

\vskip 0.3 cm

 Nous donnons trois variantes du th\'eor\`eme de Rao.

\th {Th\'eor\`eme 3.1}.  Soient $\cC$ et $\cC'$ deux familles plates de courbes
param\'etr\'ees par  $A$. Alors, $\cC$ et
$\cC'$ sont dans la m\^eme classe de biliaison si et seulement si $\cJ_\cC$ et
$\cJ_{\cC'}$  sont pseudo-isomorphes, \`a d\'ecalage pr\`es (i.e., il
existe un entier $h$ tel que $\cJ_\cC$ et
$\cJ_{\cC'}(h)$  sont pseudo-isomorphes).

\th {Th\'eor\`eme 3.2}. Soient $\cC$ et $\cC'$ deux familles plates de courbes
param\'etr\'ees par  $A$, munies de r\'esolutions de type \sN\
(resp. \sE), avec des faisceaux $\cN$, $\cN'$ (resp. $\cE$, $\cE'$). Alors, $\cC$ et
$\cC'$ sont dans la m\^eme classe de biliaison si et seulement si $\cN$ et $\cN'$
sont pseudo-isomorphes \`a d\'ecalage pr\`es (resp. si $\cE^\vee$ et $\cE'^\vee$
sont pseudo-isomorphes \`a d\'ecalage pr\`es).

\th {Th\'eor\`eme 3.3}. Soient $\cC$ et $\cC'$ deux familles plates de courbes
param\'etr\'ees par  $A$, munies de r\'esolutions de type \sN\
extraverties (resp. de type \sE\ introverties), avec des faisceaux $\cN$, $\cN'$
(resp.
$\cE$,
$\cE'$). Alors,
$\cC$ et
$\cC'$ sont dans la m\^eme classe de biliaison si et seulement si $\cN$ et $\cN'$
(resp. si $\cE$ et $\cE'$)
sont stablement isomorphes \`a d\'ecalage pr\`es.

\tarte {b) D\'emonstration : condition n\'ecessaire}

On suppose $\cC$ et $\cC'$ dans la m\^eme classe de biliaison. 

Vu les d\'efinitions
de la biliaison et des \ps, on peut supposer qu'on passe de $\cC$ \`a $\cC'$ par une
biliaison \'el\'ementaire. On a alors une surface $\cQ$ contenant $\cC$ et $\cC'$ et
un entier $h$ tels que l'on ait un isomorphisme $\cJ_{\cC/\cQ} \simeq
\cJ_{\cC'/\cQ}(h)$, de sorte que ces faisceaux sont pseudo-isomorphes. On en
d\'eduit  que $\cJ_\cC$ et $\cJ_{\cC'}(h)$ le sont aussi par 2.18.2, d'o\`u le sens
direct de 3.1. On en d\'eduit aussi, par 2.18.1, que les faisceaux $\cN$ et
$\cN'(h)$ de r\'esolutions de type \sN\ de $\cC$ et $\cC'$ sont pseudo-isomorphes,
d'o\`u 3.2 pour le type \sN\ et 3.3 pour le type \sN\ extraverti r\'esulte alors de
2.13. Pour les assertions sur les r\'esolutions de type \sE\ on lie $\cC$ (resp.
$\cC'$) \`a $\cC_1$ (resp. $\cC'_1$) et on applique ce qui pr\'ec\`ede \`a ces
courbes en utilisant 1.8 et 2.25.

\tarte {c) D\'emonstration : condition suffisante}

Notons d\'ej\`a que, quitte \`a faire deux liaisons comme ci-dessus, on peut se
contenter de montrer 3.3 (ou 3.2) dans le cas des r\'esolutions de type \sN. On
en d\'eduira alors  3.1 gr\^ace \`a 2.18.

On va montrer la proposition suivante : 

\th {Proposition 3.4}. Soient $\cC$ et $\cC'$ deux courbes de $\bP^3_A$ admettant les
r\'esolutions de type \sN\ (ou de type \sN\ d\'ecal\'ee) suivantes :
$$0 \fl \cP \fl \cN \fl \cJ_\cC \fl 0 \quad {\rm et} \quad 0 \fl \cP' \fl \cN \fl
\cJ_{\cC'}(h) \fl 0$$
o\`u $h$ est un entier et o\`u le  faisceau $\cN$ est le m\^eme pour les deux
r\'esolutions. Alors,
$\cC$ et
$\cC'$ sont dans la m\^eme classe de biliaison.

Il est clair que cette proposition implique la condition suffisante dans 3.3 ou 3.2. En
effet, si on a des r\'esolutions de type \sN\ extraverties de $\cC$ et $\cC'$ avec
$\cN$ et $\cN'$ stablement isomorphes \`a d\'ecalage pr\`es on obtient des
r\'esolutions de type \sN\ avec le m\^eme $\cN$ comme en 3.4 en ajoutant des
dissoci\'es, cf. remarque 2.17.1.

\tarte {d) D\'emonstration de 3.4}

En passant aux sections globales on a les suites exactes de modules d\'eduites des
r\'esolutions ci-dessus :
$$0 \fl P \Fl {\f} N \Fl {p} I_\cC \fl 0 \quad {\rm et}Ê\quad 0 \fl P' \Fl {\f'} N \Fl
{p'} I_{\cC'} (h) \fl 0,$$
que l'on peut r\'ecrire  sous la forme :
$$0 \fl P\oplus P' \Fl{\psi} N\oplus P' \Fl {(p,-p\f')} I_\cC \fl 0 \quad {\rm et}
\quad 0 \fl P\oplus P'
\Fl{\psi'} N \oplus P
\Fl {(p',-p'\f)}
I_{\cC'} (h) \fl 0 $$
$${\rm avec} \qquad \psi = \pmatrix {\f&\f'\cr 0& \Id_{P'} \cr} \qquad {\rm et} \qquad
\psi'=\pmatrix {\f&\f'\cr  \Id_{P}&0 \cr}$$ 
ce qui nous ram\`ene \`a montrer la proposition  suivante :

\th {Proposition 3.5}. Soient $\cC$ et $\cC'$ deux familles de courbes avec des
r\'esolutions de la forme :
$$ 0 \fl P \ \Fl {^t(\f\ \a)}\ N \oplus L \Fl{p} I_{\cC}(h) \fl 0 \quad {\rm et} \quad 
0 \fl P\ \Fl {^t(\f\ \a')}\ N \oplus L' \Fl{p'} I_{\cC'}(h') \fl 0$$
o\`u $P,L,L'$ sont libres sur $R_A$, o\`u le faisceau associ\'e \`a $N$ est
localement libre  et o\`u la fl\`eche $\f : P \fl N$ est la m\^eme pour les deux
suites. Alors, on passe de
$\cC$
\`a
$\cC'$ par un nombre fini de biliaisons \'el\'ementaires.

\dem (de 3.5) Elle proc\`ede par r\'ecurrence sur $\rg L= \rg L'$ et
n\'ecessite plusieurs lemmes. Dans tous ces lemmes les hypoth\`eses sur les modules
$P$ et $N$ seront celles de 3.5 :
$P$ est libre sur
$R_A$ et le faisceau associ\'e \`a $N$ est localement libre sur $\bP^3_A$.

 Si $s $ est une section d'un
faisceau $\cF$ sur $\bP^3_A$ on note $\ov s$ l'image de $s$ dans $\cF \T_A k(t)$. On
pose $R= R_A \T_A k(t)$.

\vskip 0.3 cm

 Rappelons, cf. 1.6, que si on a une
courbe $\cC$ contenue dans une surface  $\cQ$ d'\'equation $Q$, plate sur $A$,   un
entier
$h \geq 0$ et  un
\'el\'ement  $H \in H^0 (\bP^3,
\cO_{\bP} (h))$ tel que la multiplication par $H$,  $u_H :\cO_{\cQ}(-h)\fl
\cO_{\cQ}$, soit injective  dans
la fibre ferm\'ee,  on obtient 
 une courbe $\cC'$, avec une  biliaison
\'el\'ementaire   triviale de $\cC$ \`a $\cC'$, en prenant pour $\cO_{\cC'}$ le
conoyau de
$u_H$ et qu'alors l'id\'eal de $\cC'$ est $H I_\cC + (Q)$.  La traduction de ce
r\'esultat est le lemme suivant :

\th {Lemme 3.6}. Soit $\cC$ une famille de courbes admettant une
r\'esolution :
$$0 \fl P \Fl{\f} N\Fl{p} I_\cC \fl 0.  $$
Soient $a$ et $h$ des entiers avec $a>0$ et $h \geq 0$ et soit $Q\in H^0 (\bP^3_A,\cJ_\cC
(a))$  tel que $\ov Q \neq 0$. Alors la
famille de surfaces
$\cQ$ d\'efinie par
$Q$ est  plate sur $A$ et il existe $H\in H^0 (\bP^3_A,\cO_{\bP}
(h))=R_{A,h}$ tel que
$\ov H$ et 
$\ov Q$ soient sans facteur commun dans $R$.
   L'\'el\'ement  $H$ permet de  faire une biliaison \'el\'ementaire
triviale de hauteur $h$ sur
$\cQ$ et on obtient
 une famille de courbes $\cC'$ qui a une r\'esolution :
$$0 \fl P  \oplus R_A(-a)\Fl{\psi} N\oplus R_A(h-a)\Fl{(H p, Q)} I_{\cC'}(h) \fl
0  $$
avec $\psi=\pmatrix{\f & f\cr 0&-H\cr}$ o\`u $f$ v\'erifie $p(f)=Q$.

\dem
Puisque $\ov Q$ est non nul, $Q$ n'est pas diviseur de 0 dans $R_A$ et la famille
de surfaces d\'efinie par $Q$ est  plate sur $A$. Il est clair qu'on peut trouver un
$H$ convenable. Puisque
$\ov H$ et 
$\ov Q$ n'ont pas de  facteur
commun, la multiplication par $H$ : $\cO_{\cQ}(-h)\fl \cO_{\cQ}$
est injective  dans
la fibre ferm\'ee. On peut donc faire une biliaison \'el\'ementaire
triviale et on obtient
 une famille de courbes $\cC'$ telle que $  I_{\cC'} =H I_\cC+(Q)$.
On  \'ecrit alors $Q=p(f)$ avec $f \in N_a$ et on v\'erifie que $ I_{\cC'} $ a la
r\'esolution annonc\'ee.

\vskip 0.3 cm

Le lemme suivant est le cas particulier de 3.5 o\`u l'on suppose $L$ et $L'$ de rang
$1$ et
$\det
\f$  non nul au point ferm\'e.

\th {Lemme 3.7}. Soient $\cC$, $\cC'$ des familles de courbes
avec des
r\'esolutions de la forme :
$$ 0 \fl P \ \Fl {^t(\f\ \a)}\ N \oplus R_A(-a) \Fl{p} I_{\cC} \fl 0,$$
$$0 \fl P\ \Fl {^t(\f\ \a')}\ N \oplus R_A(-a')  \Fl{p'} I_{\cC'}(a-a')\fl 0$$
avec $a' \leq a$.
Soit  $Q \in R_A$  un d\'eterminant de $\f$ ($Q$ est un \'el\'ement homog\`ene de
degr\'e $a$ de
$R_A$, d\'efini
\`a un \'el\'ement inversible de
$A$ pr\`es). On a les propri\'et\'es suivantes :
\lign
1) Les fl\`eches
$p$ et
$p'$  sont de la forme
$p = (\s\ Q)$ et
$p'= (\s'\ Q)$. \lign 2) Si, de plus,  on a $\ov Q\neq 0$, on passe de $\cC$ \`a
$\cC'$ par un nombre fini de biliaisons
\'el\'ementaires.

\dem 1)  En vertu du th\'eor\`eme de
Hilbert-Burch les fl\`eches $p$ et
$p'$ sont donn\'ees localement par les mineurs de
$^t(\f\
\a)$ et
$^t(\f\ \a')$, donc sont de la forme annonc\'ee.

2) Supposons d'abord qu'on a $a'=a$. Puisque $\ov Q$ est non nul, la
famille de surfaces
$\cQ$ d\'efinie par
$Q$ est plate sur $A$.
Les id\'eaux
$I_{\cC/\cQ}$ et
$I_{\cC'/\cQ}$ ont alors tous deux une r\'esolution de la forme
$$0 \fl P \Fl{\f}  N \fl I_{\cC/\cQ} \fl 0,$$  donc sont tous deux
isomorphes \`a
$\Coker \f$ et on a un isomorphisme
$u : I_{\cC/\cQ}\fl I_{\cC'/\cQ}$  c'est-\`a-dire encore une
biliaison \'el\'ementaire de hauteur $0$ qui passe de
$\cC$ \`a $\cC'$.

Si $a'$ est $<a$ on prend $H \in H^0(\bP^3_A, \cO_\bP (a-a'))$ tel que $\ov H$ et $\ov Q$
soient sans facteur commun et on effectue la biliaison \'el\'ementaire triviale
de hauteur $h=a-a'$ d\'efinie par
$Q$ et
$H$ (cf. 3.6). On obtient une courbe $\cC_1$ avec la r\'esolution donn\'ee par 3.6 :
$$ 0 \fl P \oplus R_A(-a) \Fl {\psi} N \oplus R_A(-a) \oplus R_A(-a') \fl I_{\cC_1} (a-a')
\fl 0$$
avec 
$\psi = \pmatrix { \f & 0 \cr \a & 1 \cr 0 & -H \cr}$. 
On v\'erifie  qu'on peut simplifier le terme $R_A (-a)$ pour obtenir  la
r\'esolution :
$$0 \fl P \Fl {^t (\f, H \a)} N \oplus R_A(-a') \Fl {(H \s, Q)} I_{\cC_1} (a-a') \fl 0.$$
La conclusion r\'esulte alors du cas $a=a'$ appliqu\'e \`a $\cC'$ et $\cC_1$.

\vskip 0.5 cm

Nous abordons maintenant le cas $\ov Q=0$.

\th {Lemme 3.8}. Soit $\cC$ une  famille de courbes avec une
r\'esolution de la forme :
$$ 0 \fl P \ \Fl {^t(\f\ \a)}\ N \oplus R_A(-a) \Fl{(\s\ Q)} I_{\cC} \fl 0.$$ On
suppose qu'on a  $\ov Q= 0$.\lign Il existe
$b\in \bN$ tel que, pour $(s,H,H',\beta)$ g\'en\'eral dans $N_b\times
R^2_{A,b-a}\times \Hom_{R_A} (P,R_A(-a))$,
les conditions suivantes sont r\'ealis\'ees :\lign
(i)  les images de $ \det \pmatrix {\f &s\cr \a &H\cr}$ et $\det \pmatrix {\f &s\cr \beta
&H'\cr}$ dans $R$ sont non nulles, \lign
(ii) $\ov{ H'}$ et $\ov {\s (s)+Q H}$ sont sans facteur commun, \lign
(iii) il existe une famille plate $\cC'$ de courbes admettant  la r\'esolution :
 $$0 \fl  P\oplus
R_A(-b)\Fl {\psi} N
\oplus R_A(-a)^2 \fl I_{\cC'}(b-a) \fl 0$$
 avec  $\psi =\pmatrix {\f &s\cr \a &H\cr \beta &H'\cr}$ et on passe de $\cC$ \`a $\cC'$
par un nombre fini de biliaisons \'el\'ementaires.

\dem Chaque condition \'etant ``ouverte modulo $m_A$'' (pour {\sl (i)} et {\sl (ii)}
c'est imm\'ediat, pour {\sl (iii)},
cf. par exemple [MDP1] II 1.2), il suffit de v\'erifier s\'epar\'ement  qu'elle est non
vide.

Soient $C$ la fibre de $\cC$ au-dessus du point ferm\'e et $\cN$ le faisceau
associ\'e \`a $N$. Remarquons d'abord que l'homomor\-phis\-me
$\ov\cN\to\cJ_C$ induit par $\s$ n'est pas nul puisque $C$ est une courbe.
En vertu du th\'eor\`eme de changement de base il existe $b>a$ tel que $H^0\cN(b)=N_b$, que
l'homomorphisme $H^0\cN(b)\to H^0\ov\cN(b)$ soit surjectif et que $H^0\ov\cN(b)\to
H^0\cJ_C(b)$ soit non nul, donc il existe $s\in N_b$ tel que $\ov\s
(\ov s) \neq 0$.
Le d\'eterminant de $ \pmatrix {\f &s\cr \a &H\cr}$ est \'egal \`a $\s
(s)+QH= F$ et n'est pas nul dans $R$ puisque
$\ov Q= 0$. En prenant aussi $\beta = \a$ et $H'=H$, on voit que  la condition
{\sl (i)} n'est pas vide.

Soient $s$ et $H$ v\'erifiant {\sl (i)}.  On
choisit $H'$ dans $R_{A, b-a}$ tel que $\ov F$ et $\ov H'$ soient sans facteur
commun (cf. 3.6) ce qui donne {\sl (ii)} et on fait avec $H'$ une biliaison
\'el\'ementaire triviale de hauteur
$b-a$ sur la famille de surfaces d\'efinie par $F$. Le lemme 3.6	 donne une courbe
$\cC'_1$ admettant la r\'esolution annonc\'ee en {\sl (iii)} avec $\beta=0$. Pour
$\beta$ g\'en\'eral on a encore une courbe $\cC'$. En vertu de 3.7, comme $ \det \pmatrix
{\f &s\cr
\a &H\cr}$ est non nul, on passe de $\cC'_1$ \`a
$\cC'$ (donc aussi de $\cC$ \`a $\cC'$)  par un nombre fini de biliaisons
\'el\'ementaires.
\vskip 0.3 cm

Nous pouvons maintenant prouver 3.5 dans le cas $L= L' = R_A(-a)$.

\th {Lemme 3.9}.  Soient $\cC$ et $\cC'$ deux familles de courbes
avec des
r\'esolutions de la forme :
$$ 0 \fl P \ \Fl {^t(\f\ \a)}\ N \oplus R_A(-a) \Fl{(\s\  Q)} I_{\cC} \fl 0$$
$$0 \fl P\ \Fl {^t(\f\ \a')}\ N \oplus R_A(-a)  \Fl{(\s'\ Q)} I_{\cC'} \fl 0.$$
Alors on passe
de $\cC$ \`a $\cC'$ par un nombre fini de  biliaisons \'el\'ementaires.

\dem  Si $\ov Q$ est non nul on applique 3.7. Sinon, il existe
$b\in \bN$ et $(s,H,H',\beta)$  dans $N_b\times R^2_{A,b-a}\times \Hom_{R_A}
(P,R_A(-a))$ qui v\'erifient
les conditions du lemme 3.8 pour $\cC$ et pour $\cC'$. 

Soit $\psi$ (resp.
$\psi'$) : $ P\oplus R_A(-b)\to
N\oplus R_A(-a)^2$ donn\'e par la matrice
 $$\pmatrix {\f &s\cr \a &H\cr \beta &H'\cr} \qquad \qquad {\rm (resp.} \; \pmatrix
{\f &s\cr
\a' &H\cr \beta &H'\cr}).$$ Alors $\psi$
(resp. $\psi'$) d\'efinit une famille de courbes $\cC_1$ (resp. $\cC'_1$).
Puisque le d\'eterminant de
$\pmatrix {\f &s\cr
\beta &H'\cr}$ n'est pas nul dans $R$, on passe
de $\cC_1$ \`a $\cC'_1$  par une biliaison \'el\'ementaire en vertu de 3.7. 
Mais, d'apr\`es 3.8, on passe aussi de $\cC$ \`a
$\cC_1$ (resp. de $\cC'$ \`a
$\cC'_1$) par une biliaison
\'el\'ementaire, d'o\`u le r\'esultat.

\vskip 0.3 cm
Il ne reste plus qu'\`a terminer la d\'emonstration de 3.5.

\vskip 0.3 cm

On pose $L = R_A(-a_1) \oplus \cdots R_A(-a_r)$ et $L' = R_A(-b_1) \oplus \cdots
R_A(-b_r)$ avec $a_1 \leq a_2 \leq \cdots \leq a_r$ et $b_1 \leq b_2 \leq \cdots \leq
b_r$ et on peut supposer $a_r \leq b_r$. Soient
$Q'_1,
\cdots, Q'_r$ les images des
$R_A (-b_j)$ dans
$ I_{\cC'}(h')$. On raisonne par r\'ecurrence sur $r$. Le cas
$r=0$ est
\'evident. Supposons  l'assertion \'etablie pour $r-1$.

{\it 1) Premier cas :} On suppose  qu'il existe $i,j$ tels que $a_i=b_j$.

Soient
$\a_i$ et $\a'_j$ les lignes  d'indice $i,j$ dans $\a$ et $\a'$.  Il
existe  un homomorphisme $w : P \fl R_A(-a_i)= R_A(-b_j)$ tel que si on remplace
$\a_i$ (resp.
$\a'_j$) par $w$ dans $\a$ (resp. $\a'$) le conoyau du morphisme est encore
l'id\'eal
tordu d'une famille de courbes  $\G$ (resp. $\G'$). En effet l'ensemble des
$w$ convenables est un ouvert non vide de $\Hom (P, R_A(-a_i))\T_A k$.
Alors, le lemme 3.9 montre que
$\cC$ et
$\G$ (resp.
$\cC'$ et
$\G'$) sont dans la m\^eme classe de biliaison et l'hypoth\`ese de
r\'ecurrence montre
qu'il en est de m\^eme de
$\G$ et
$\G'$, ce qui permet de conclure.

\vskip 0.3 cm

{\it 2) Deuxi\`eme cas :} On suppose $a_r < b_r$ et l'un des $\ov {Q'_j}$
non nul.

Comme on a $b_j \leq b_r$, quitte \`a faire un changement de base dans $L'$ (ce qui ne
change  pas
$\f$),  on peut
supposer $\ov {Q'_r}$ non nul.

On effectue  une biliaison \'el\'ementaire triviale sur la surface
(plate) d\'efinie par
$Q'_r$, de hauteur $b_r-a_r$ et on obtient une courbe $\cC''$. En utilisant 3.6 et un
raisonnement de simplification comme celui utilis\'e dans la preuve de 3.7, on voit que
$\cC''$ a une r\'esolution de la m\^eme forme que
$\cC'$, avec le m\^eme
$\f$, mais o\`u
$b_r$ est
remplac\'e par $a_r$ et on est ramen\'e au premier cas (attention $a_r$ n'est plus
n\'ecessairement le plus grand des $b_j$).

\vskip 0.3 cm

{\it 3) Troisi\`eme cas : }  On suppose $a_r < b_r$ et les  $\ov {Q'_j}$
tous nuls.

Soient $\ov\cC$ et $\ov \cC'$ les fibres ferm\'ees des familles $\cC$ et
$\cC'$.
On a  la suite exacte $\ov {\cP }\Fl {\ov{\f}} \ov {\cN } \fl
\cJ_{\ov\cC'}(h')\fl 0$ et, avec la
r\'esolution de $\cJ_{\ov\cC}(h)$ on en d\'eduit qu'on a une surjection
$\cJ_{\ov\cC}(h) \fl \cJ_{\ov\cC'}(h')$ ce qui n'est possible que si
$\ov\cC=\ov \cC'$ et $h=h'$. Mais
alors, si on regarde les r\'esolutions de
$\ov\cC$ et
$\ov
\cC'$ obtenues par r\'eduction modulo $m_A$, l'\'egalit\'e des cohomologies de
$\ov\cC$ et $\ov \cC'$  montre  l'\'egalit\'e des
$a_i$ et $b_i$, ce qui est absurde.

\rema {3.10}  La preuve de la condition
n\'ecessaire est encore valide lorsque le corps r\'esiduel $k(t)$ de $A$ est fini. En
revanche pour la condition suffisante, les assertions d'existence de la preuve de 3.4 utilisent
l'hypoth\`ese que $k(t)$ est infini. Si $k(t)$ est fini, la conclusion de 3.4 est vraie sur un
anneau
$B$ local fini et \'etale sur $A$. 

\vskip 0.3 cm

\tarte {e) Liaison impaire}

Les r\'esultats pr\'ec\'edents permettent aussi de donner un th\'eor\`eme dans le cas
de la liaison impaire :

\th {Th\'eor\`eme 3.11}. On suppose $A$ local. Soient $\cC$ et $\cC'$ deux familles
plates de courbes param\'etr\'ees par  $A$, munies de r\'esolutions de type \sN\
(resp. \sE), avec des faisceaux $\cN$, $\cN'$ (resp. $\cE$, $\cE'$). Alors, $\cC$ et
$\cC'$ sont li\'ees par un nombre impair de liaisons \'el\'ementaires si et
seulement si
$\cN$ et
$\cE'^\vee$ sont pseudo-isomorphes  \`a d\'ecalage pr\`es ou encore si et
seulement si $\cN'$ et
$\cE^\vee$ sont pseudo-isomorphes  \`a d\'ecalage pr\`es.

\dem Supposons d'abord $\cC$ et $\cC'$ reli\'ees par une cha\^\i ne de courbes $\cC_0
= \cC, \cC_1, \cdots, \cC_{2n}, \cC_{2n+1} = \cC'$, dans laquelle on passe de $\cC_i$
\`a $\cC_{i+1}$ par une liaison \'el\'ementaire. Alors, on sait (cf. 2.24) que
$\cC_{2n}$ admet une r\'esolution de type \sE\ dont le faisceau $\cE_{2n}$ est 
$\cN'^\vee (-s-t)$. Comme $\cC$ et $\cC_{2n}$ sont bili\'ees en vertu de 1.8, il
r\'esulte de 3.2 que
$\cE^\vee$ et $\cN' (-s-t)$ sont pseudo-isomorphes \`a d\'ecalage pr\`es. Pour l'autre
assertion on applique 2.24 \`a $\cC$ et $\cC_1$.

Supposons maintenant, par exemple, $\cN$ et
$\cE'^\vee$ sont pseudo-isomorphes  \`a d\'ecalage pr\`es. On fait une
liaison qui passe de $\cC$ \`a $\cC_1$ (cf. 1.4). Alors, $\cC_1$ a une r\'esolution de
type \sE\ avec pour $\cE_1$ un d\'ecal\'e de $\cN^\vee$ (cf. 2.24) et donc $\cN$
pseudo-isomorphe \`a $\cE'^\vee$ (cf. 2.3). Il r\'esulte alors de 3.2 que $\cC_1$ et
$\cC'$ sont bili\'ees, donc que l'on passe de l'une \`a l'autre par un nombre pair de
liaisons
\'el\'ementaires (cf. 1.8), donc de $\cC$ \`a $\cC'$ par un nombre impair de
liaisons \'el\'ementaires.

\vskip 0.3 cm

\tarte {f) Cas de la dimension quelconque}

La d\'efinition des pseudo-isomorphismes se g\'en\'eralise comme suit dans le cas de
la dimension $n \geq 3$ :

\th {D\'efinition 3.12}. Soient $\cN$ et $\cN'$  des faisceaux coh\'erents sur $\bP^n_A$
et plats sur $A$  et soit 
$f$ un morphisme de 
$\cN$ dans $\cN'$. On dit que $f$ est un {  pseudo-isomorphisme} s'il induit
:\lign 0) un isomorphisme de foncteurs $H^0(\cN(n)\T_A \tas) \fl H^0(\cN'(n)\T_A \tas)$ pour
tout $n \ll 0$,  \lign
1) un isomorphisme de foncteurs $H^i_*(\cN\T_A \tas) \fl
H^i_*(\cN'\T_A \tas)$ pour 
$1 \leq i \leq n-2$ et
\lign 2) un monomorphisme de foncteurs $H^{n-1}_*(\cN\T_A \tas) \fl H^{n-1}_*(\cN'\T_A
\tas)$.\lign Deux faisceaux coh\'erents sur $\bP^n_A$ et plats sur $A$ seront dits
pseudo-isomorphes s'il existe une cha\^\i ne de \ps\  qui
les joint :
$$ \cN= \cN_0 \fl \cN_1 \lf \cN_2 \fl \cN_3 \lf \cdots \fl \cN_{2p-1} \lf \cN_{2p} =
\cN'.$$

On g\'en\'eralise sans difficult\'e l'ensemble des d\'efinitions et  des r\'esultats
ci-dessus  et on obtient le th\'eor\`eme suivant :

\th {Th\'eor\`eme 3.13}. Soient $\cC$ et $\cC'$ deux familles plates de sous-sch\'emas
de $\bP^n$ ($n \geq 3$),
de pure codimension $2$, sans composantes immerg\'ees, 
param\'etr\'ees par l'anneau local
$A$, munies de r\'esolutions de type
\sN\ (resp. \sE), avec des faisceaux $\cN$, $\cN'$ (resp. $\cE$, $\cE'$). Alors, $\cC$ et
$\cC'$ sont dans la m\^eme classe de biliaison si et seulement si $\cN$ et $\cN'$
sont pseudo-isomorphes,  \`a d\'ecalage pr\`es (resp. si $\cE^\vee$ et $\cE'^\vee $
sont pseudo-isomorphes,  \`a d\'ecalage pr\`es).

\vskip 1 cm

\titre {R\'ef\'erences bibliographiques}

[AG] Hartshorne R., Algebraic geometry, Graduate texts in Mathematics 52, Springer
Verlag, 1977.

[BPS] B\u anic\u a C., Putinar M., Schumacher G., Variation der globalen Ext in
Deformationen kompakter komplexer R\" aume, Math. Ann. 250, 1980, 135-155.

[EGA] Grothendieck A. et Dieudonn\'e J., \'El\'ements de g\'eom\'etrie alg\'ebrique,
Publ. math. IHES, 1960-1967.

[H] Hartshorne R., Coherent functors, \`a para\^\i tre, Advances in Math.

[HMDP2] Hartshorne R., Martin-Deschamps M. et Perrin D., Construction de familles
minimales de courbes gauches, en pr\'eparation.

[HMDP3] Hartshorne R., Martin-Deschamps M. et Perrin D., Triades et familles
de courbes gauches, en pr\'eparation.

[Ho] Horrocks G., Vector bundles on the punctured spectrum of a local ring, Proc. London
Math. Soc. (3) 14, 1964, 689-713.

[JS] de Jong T., van Straten D., Deformations of the normalization of hypersurfaces 
(Appendix), Math. Ann. 288, 1990, 527-547.

 [K] Kleppe J.-O., Liaison of families of
subschemes in
$\bP^n$  in Algebraic curves and projective geometry, (128-173), Proc. Trento, 1988,
Lecture Notes in Math. 1389, Springer Verlag.  

[MDP1]  Martin-Deschamps M. et  Perrin D., Sur la classification des
courbes gauches, Ast\'erisque, Vol. 184-185, 1990.

[Mi] Migliore J., An introduction to deficiency modules and liaison theory for
subschemes of projective space, Lecture notes series of research institute of math.,
Seoul National University, 1994.

[N] Nollet S., Even linkage classes, Trans. Amer. Math. Soc., 348-3, 1996, 1137-1162.

[PS] Peskine C. et Szpiro L., Liaison des vari\'et\'es alg\'ebriques, Invent. Math.,
26, 1974, 271-302.

[R1] Rao A.P., Liaison among curves in $\bP^3$, Invent. Math., 50, 1979,
205-217.

[R2] Rao A.P., Liaison equivalence classes, Math. Ann. 258, 1981, 169-173.

[S] Schlesinger E., Ext and base change, unpublished manuscript, Berkeley, 1993.

[V] Verdier J.-L., Cat\'egories d\'eriv\'ees, \'etat 0, {\it in} SGA $4 {1 \over 2}$,
Lecture Notes in Math. 569, Springer Verlag, 1977.

\bye